\documentclass[a4,12pt]{article}

\usepackage{amssymb,amsmath}

\newtheorem{defn}{Definition}[section]
\newtheorem{proposition}[defn]{Proposition}
\newtheorem{corollary}[defn]{Corollary}
\newtheorem{rem}[defn]{Remark}
\newtheorem{exm}[defn]{Example}
\newtheorem{lemma}[defn]{Lemma}
\newtheorem{theorem}[defn]{Theorem}
\newtheorem{xproof}{{\it Proof. }}

\newenvironment{definition}{\begin{defn}\em}{\end{defn}}
\newenvironment{remark}{\begin{rem}\em}{\end{rem}}
\newenvironment{example}{\begin{exm}\em}{\end{exm}}

\newenvironment{proof}{\begin{xproof}\em}{\end{xproof}}

\def\qed{\hspace{0.3cm}{\rule{1ex}{2ex}}}

\newcommand\cf{cf.}
\newcommand\opens[1]{\tau_{#1}}
\newcommand\CC{\mathbb{C}}
\newcommand\orth{\mathrm{orth}}
\newcommand\ann{\mathrm{ann}}
\newcommand\suplattices{{\bf SL}}
\newcommand\V{\bigvee}
\newcommand\Q{{\cal Q}}
\newcommand\ie{i.e.}
\newcommand\eg{e.g.}
\newcommand\Max{{\rm Max}\,}
\newcommand\then{\!\cdot\!}
\newcommand\pwset[1]{2^{#1}}
\newcommand\st{\mid}
\newcommand\opp[1]{{#1}^\textrm{op}}
\newcommand\rs{{\rm R}}
\newcommand\ls{{\rm L}}
\newcommand\rorth[1]{{{#1}^{\perp}}}
\newcommand\lorth[1]{{^{\perp}{#1}}}
\newcommand\upsegment{{\uparrow}}
\newcommand\downsegment{{\downarrow}}
\newcommand\sub{\mathrm{Sub}}
\newcommand\csub{\overline{\mathrm{Sub}}}
\newcommand\form[2]{\langle #1,#2\rangle}
\newcommand\B{\mathcal{B}}
\newcommand\proj{\mathcal{P}}

\begin{document}

\title{Sup-lattice 2-forms and quantales\thanks{Research
partially supported by FCT and FEDER
via the research center CLC of IST and through grant
POCTI/1999/MAT/33018.}}
\author{Pedro Resende
\vspace*{1mm}\\
\small\it Departamento de Matem{\'a}tica,
Instituto Superior T{\'e}cnico,
\vspace*{-1mm}\\
\small\it Av. Rovisco Pais, 1049-001 Lisboa, Portugal
\vspace*{1mm}\\
\small\rm E-mail: pmr@math.ist.utl.pt}

\date{~}

\maketitle

\begin{abstract}
A 2-form between two sup-lattices $L$ and $R$ is defined to be a 
sup-lattice bimorphism $L\times R\rightarrow 2$. Such 2-forms are 
equivalent to Galois connections, and we study them and their 
relation to quantales, involutive quantales and quantale modules. 
As examples we describe applications to C*-algebras.
\vspace{0.1cm}
\vspace{0.1cm}
\end{abstract}

\section{Introduction}

Let $L$ and $R$ be sup-lattices. A Galois connection between $L$ 
and $R$ is a pair of antitone maps $\rorth{(-)}:L\rightarrow R$
and $\lorth{(-)}:R\rightarrow L$ such that
$x\le \lorth{(\rorth x)}$ and $y\le\rorth{(\lorth y)}$ for all $x\in 
L$ and $y\in R$.
In fact all the information present in the Galois connection is
already
available in each of the maps $\rorth{(-)}$ and $\lorth{(-)}$, 
due to completeness of the lattices, 
or equivalently in the map
$\varphi:L\times R\rightarrow 2$ given by
\[\varphi(x,y)=0\iff 
x\le\lorth y\ \ \ (\iff y\le\rorth x)\;,\]
which is a 
bimorphism of sup-lattices. In this paper we 
study such bimorphisms, and call them (sup-lattice) 
2-forms.
The purpose is to provide a useful framework within which
to study various 
aspects of quantales and their modules, including involutive 
quantales and their applications to C*-algebras.

We study two notions of map between 2-forms: the 
orthomorphisms, in \S\ref{sec:2forms}, which are analogous to
isometries, and the 
continuous maps, in \S\ref{sec:quantalesandtwoforms},
so-named because they generalize the
continuous maps of topological spaces. In particular,
the set of continuous 
endomaps of a 2-form has the structure of a quantale, we show that
under mild 
restrictions the 2-form can be recovered from it, and we 
obtain generalizations of well known facts~\cite{MulvPell92}
concerning the right and 
left sides of quantales of sup-lattice 
endomorphisms, and also concerning involutive quantales. In 
\S\ref{sec:principalqumods} we deal with principal quantale modules
(\ie, modules with a 
single generator), and in \S\ref{sec:modsand2forms} we relate them
to 2-forms. 
Finally, in \S\ref{sec:invmods} we address the
particular case of symmetric 
2-forms and involutive quantales, and we discuss applications to 
C*-algebras.

We are indebted to the work of Mulvey and 
Pelletier~\cite{MulvPell01}, which was one of the main sources
of inspiration for 
our paper. They implicitly use parts of the theory
of 2-forms, 
and this is reflected in the fact that we obtain, in
\S\ref{sec:invmods}, 
a much shorter proof of one of their main 
theorems~\cite[Th.\ 9.1]{MulvPell01}, which concerns the relation 
between quantales and C*-algebras.
We hope in this way to bring out more 
explicitly some of the principles that lie behind that relation.

\section{Background}
\label{sec:background}

In this section we present some basic facts, terminology and
notation
concerning sup-lattices, quantales and quantale modules,
however without 
attempting to be complete. Further basic reading about
sup-lattices and
quantales can be found in the first chapters of the book by
Rosenthal~\cite{Rose90},
and further references will be
cited throughout this section.

By a \emph{sup-lattice} is meant a partially ordered set $S$ each of
whose subsets $X\subseteq S$ has a join (supremum) $\V X$ in $S$
(hence, a sup-lattice is a complete lattice). By a \emph{homomorphism}
of sup-lattices $f:S\rightarrow T$ is meant a map that preserves
arbitrary 
joins: $f(\V X)=\V\{f(x)\st x\in X\}$, for all $X\subseteq S$.
The greatest element $\V S$ of a sup-lattice $S$ (the top) is
denoted by
$1_{S}$, or
$1$, and the least element $\V\emptyset$ (the bottom) by
$0_{S}$, or $0$. The two-element
sup-lattice $\{0,1\}$ is denoted by $2$. The order-dual of a 
sup-lattice $S$, \ie, $S$ with the order reversed, is denoted
by $\opp S$.
A homomorphism of
sup-lattices
$f:S\rightarrow T$ is said to be \emph{strong} if
$f(1)=1$, and \emph{dense} if the condition $f(x)=0$ implies $x=0$ for all
$x\in S$.

Any sup-lattice homomorphism $f:S\rightarrow T$ has a right 
adjoint $f_{*}:T\rightarrow S$, which preserves all the meets
(infima) in $T$ and is defined by
\[f_{*}(y)=\V\{x\in S\st f(x)\le y\}\;.\] Equivalently,
$f_{*}$ is the unique 
monotone map that 
satisfies the condition
\[f(x)\le y\iff x\le f_{*}(y)\] for all $x\in S$ 
and $y\in T$.

If $S$ is a sup-lattice, and $j:S\rightarrow S$ is
a closure operator 
on $S$, the set of fixed-points of $j$,
$S_{j}=\{x\in S\st j(x)=x\}$,
is a sup-lattice whose joins are given by
$\V^j X=j(\V X)$, and the map $j:S\rightarrow S_{j}$
that sends each 
$x\in S$ to $j(x)$ is a surjective sup-lattice
homomorphism. Any 
quotient of a sup-lattice arises like this,
up to isomorphism, for if
$f:S\rightarrow T$ is a surjective sup-lattice
homomorphism then
$j=f_{*}\circ f$ is a closure operator on $S$,
and $T\cong 
S_{j}$~\cite{JoyaTier84,Rose90}.

The category $\suplattices$ of sup-lattices is 
monoidal~\cite{JoyaTier84},
and a 
semigroup in it is a \emph{quantale}, \emph{unital}
if the 
semigroup is a monoid, \emph{involutive} if the
semigroup has an 
involution. A left (resp.\ right) \emph{module} over a quantale $Q$ is 
a left (resp.\ right) action in
$\suplattices$.
The multiplication of two elements $a$ and $b$ in a quantale $Q$
is 
denoted by $a\then b$; if the quantale is unital, its
multiplicative unit is denoted by $e_Q$, or simply $e$;
if the quantale is involutive, the involution assigns
to each $a\in Q$ an element that is denoted by $a^*$.
The action of an element $a\in Q$
on $x\in M$, where $M$ is a left $Q$-module, is denoted
by $ax$ (or 
$xa$ for a right $Q$-module), and the module is
\emph{unital} if $ex=x$ for all $x\in M$
(resp.\ $xe=x$ for a right module). An element $a$
of a quantale is 
\emph{left-sided} (resp.\ \emph{right-sided}) if
$1\then a\le a$ (resp.\ $a\then 1\le a$). An element
which is both 
left- and right-sided is \emph{two-sided}.
The set
of left-sided elements of a 
quantale $Q$ is denoted by $\ls(Q)$ (this is a right
$Q$-module under 
multiplication), and the set (a left $Q$-module)
of right-sided 
elements is denoted by $\rs(Q)$. A \emph{factor} is a
quantale $Q$ whose set of 
two-sided elements is $\{0,1\}$.

For any sup-lattice $S$ the set $\Q(S)$ of
sup-lattice endomorphisms
of $S$ is a unital quantale under the pointwise
ordering, with 
multiplication given by composition,
$f\then g=g\circ f$, and we
have $\ls(\Q(S))\cong S$ and 
$\rs(\Q(S))\cong\opp S$~\cite{MulvPell92}. Explicitly,
for a unital quantale $Q$ we have $\ls(Q)=1\then Q$ and
$\rs(Q)=Q\then 
1$, and thus, for $\Q(S)$, a left-sided element
is the same as a ``constant'' 
map for some $s\in S$,
\[c_{s}(x)=\left\{\begin{array}{ll}
s&\mbox{if }x\neq 0\;,\\
0&\mbox{if }x=0\;,
\end{array}\right.\]
and a right-sided element is an annihilator of some
$s\in S$,
\[a_{s}(x)=\left\{\begin{array}{ll}
1&\mbox{if }x\not\le s\;,\\
0&\mbox{if }x\le s\;.
\end{array}\right.\]
It also follows from this that $\Q(S)$ is a factor.

Another example of unital
quantale, for any monoid $M$, is the powerset $\pwset M$
under pointwise multiplication:
\[X\then Y=\{xy\st x\in X,\ y\in Y\}\:.\]
This construction is universal in the sense that for any unital
quantale $Q$ and any homomorphism of monoids $h:M\rightarrow Q$
there is a unique homomorphism of unital quantales
$\bar h:\pwset M\rightarrow Q$ such that
$\bar h(\{x\})= h(x)$ for all $x\in M$. Hence, any unital quantale
is a quotient of one of the form $\pwset M$, for some monoid $M$.

Quotients of quantales and modules can be described in terms of
closure operators with additional properties: a
\emph{quantic nucleus}~\cite{Rose90}, or simply a \emph{nucleus},
on a unital quantale $Q$ is a closure operator $j$ on $Q$ such that
for all $a,b\in Q$ we have $j(a)\then j(b)\le j(a\then b)$; and a
\emph{nucleus} on a left $Q$-module $M$ is
a closure operator $k$ on $M$ such that for all $a\in Q$ and 
$x\in M$ we have $ak(x)\le k(ax)$ (see~\cite{Pase02} 
or~\cite[\S 2.5]{Rese02}). Given nuclei $j$ and $k$ as above, $Q_j$
is a unital quantale with multiplication $(a,b)\mapsto a\ast b$
defined by $a\ast b=j(a\then b)$, and $M_k$ is a left $Q$-module
with action $(a,x)\mapsto a\bullet x$ defined by $a\bullet 
x=k(ax)$. The surjective maps $j:Q\rightarrow Q_j$ and
$k:M\rightarrow M_k$ are respectively a homomorphism of unital
quantales and a homomorphism of left $Q$-modules, and any quotient
of quantales or of left $Q$-modules arises like this, up to 
isomorphism. If $j$ and $k$ are further related by the condition
$j(a) x\le k(a x)$, for all $a\in Q$ and $x\in M$, then $M_k$ is
also a left $Q_j$-module.

If $R$ is a ring with unit then the map that sends 
each subset of $R$ to the additive subgroup it generates is a
nucleus on the unital quantale $\pwset R$, and thus
$\sub(R)$, the set of additive subgroups of $R$, is a unital
quantale with 
multiplication defined by
\[a\then b=\{r_{1}s_{1}+\cdots+r_{n}s_{n}\st 
r_{i}\in a,\ s_{i}\in b\}\;.\] The left-sided elements of
$\sub(R)$ are then the left 
ideals of $R$.

More generally, in the case of a unital $k$-algebra $A$ with
$k$ an arbitrary commutative ring, the set $\sub_k(A)$ of all
the $k$-submodules of $A$ is a unital quantale, and if $A$ is
a topological $k$-algebra then the set $\csub_k(A)$
of all the closed $k$-submodules of $A$ is a unital quantale with 
multiplication defined by
\[a\then b=\overline{\{r_{1}s_{1}+\cdots+r_{n}s_{n}\st 
r_{i}\in a,\ s_{i}\in b\}}\;,\]
where $\overline {(-)}$ denotes 
topological closure.

We can obtain examples of modules in a similar way. If $R$ is a ring
and $M$ is a left $R$-module then the set $\sub(M)$ of additive
subgroups of $M$ is a left module over 
$\sub(R)$, with action defined by
\[ax=\{r_{1}m_{1}+\cdots+r_{n}m_{n}\st 
r_{i}\in a,\ m_{i}\in x\}\;.\]
A similar expression, but including closure for the 
topology, gives us a 
left 
$\csub_k(A)$-module $\csub_k(M)$, consisting of all the closed
$k$-submodules of $M$, from any topological left $A$-module $M$
over a topological $k$-algebra $A$.

If $A$ is a complex C*-algebra with unit, the unital quantale
$\csub_\CC(A)$
is involutive, with involution obtained pointwise from the
involution
of $A$. This involutive quantale is denoted by $\Max A$
in~\cite{KPRR,Mulv89,MulvPell01,MulvPell02}, where it plays the role of
the
``noncommutative maximal spectrum" of $A$. If $H$ is a Hilbert space,
its norm-closed linear subspaces can be identified with the
projections
on $H$, and we denote the sup-lattice $\csub_\CC(H)$ by $\proj(H)$.
Any C*-algebra representation $\pi:A\rightarrow\B(H)$ of $A$ on $H$
makes $H$ a topological left $A$-module, thus making $\proj(H)$ a
left $\Max A$-module.

Let $L$, $R$, and $M$ be sup-lattices, and
$(-)*(-):L\times R\rightarrow M$ a sup-lattice \emph{bimorphism},
\ie, a map 
that preserves joins in each variable (\eg, the multiplication
$Q\times 
Q\rightarrow Q$ of a quantale $Q$, or the action
$Q\times M\rightarrow M$
of $Q$ on a left module $M$):
\begin{eqnarray*}
\left(\V X\right)*y &=& \V\{x*y\st x\in X\}\;;\\
x*\left(\V Y\right) &=& \V\{x*y\st y\in Y\}\;.
\end{eqnarray*}
We will consistently use the following 
notation for the residuations associated to $*$ (\ie, the right
adjoints to the homomorphisms $(-)*y$ and $x*(-)$), for each
$x\in L$, $y\in R$, and $z\in M$:
\begin{eqnarray*}
    z/y&=&\V\{x\in L\st x*y\le z\}\;,\\
    x\backslash z&=&\V\{y\in R\st x*y\le z\}\;.
\end{eqnarray*}
Also, we define the following annihilators:
$\ann(x)=x\backslash 0$, and $\ann(y)=0/y$.
Hence, we have
\[\begin{array}{c}
y\le x\backslash z\iff x*y\le z \iff x\le z/y\;,\\
y\le\ann(x)\iff x*y=0\iff x\le\ann(y)\;,
\end{array}\]
and also the following (in)equalities:
$(z/y)*y\le z$, $\ann(y)*y=0$, $x*(x\backslash z)\le z$,
$x*\ann(x)=0$,
$x\le (x*y)/y$, $y\le x\backslash (x*y)$,
$((x*y)/y)*y=x*y$, $x*(x\backslash (x*y))=x*y$.

\section{2-forms and orthomorphisms}
\label{sec:2forms}

Let $S$ be a sup-lattice. Since its dual, $\opp S$, is order
isomorphic
to $\hom(S,2)$~\cite{JoyaTier84}, any Galois
connection between two sup-lattices $L$ and
$R$ is uniquely determined by a sup-lattice homomorphism
$L\rightarrow\hom(R,2)$, which in turn is equivalent to a
sup-lattice
bimorphism $L\times R\rightarrow 2$ (because we have an order
isomorphism $\hom(L\otimes R,M)\cong\hom(L,\hom(R,M))$ for any
sup-lattices $L,R,M$~\cite{JoyaTier84}).
Such bimorphisms are analogous to bilinear forms
on ring modules, and provide us with a convenient alternative
language
for describing Galois connections.

\begin{definition}
    Let $L$ and $R$ be sup-lattices. A map $\varphi:L\times
    R\rightarrow 2$ that preserves arbitrary joins in each variable
    is called a \emph{2-form} between $L$ and $R$, and we usually
    write
    $\form x y$ or ${\form x y}_{\varphi}$
    instead of $\varphi(x,y)$. Two elements $x\in L$ and $y\in R$
    are
    \emph{orthogonal} if $\form x y=0$, in which case we
    write $x\perp y$. The form is
    \emph{dense on the right} if $1\perp y$ implies $y=0$,
    \emph{dense on the left} if $x\perp 1$ implies $x=0$,
    and \emph{dense} if it is both dense on the right and on the
    left.
    The form
    is \emph{faithful on the right} if $x=y$ whenever
    $\form z x=\form
    z y$ for all $z\in L$, \emph{faithful on the left} if $x=y$
    whenever $\form x z=\form y z$ for all $z\in R$, and
    \emph{faithful}, or \emph{non-singular}, or a \emph{duality},
    if it is both faithful on the right and on the left.
    A 2-form $\varphi:S\times S\rightarrow 2$ is \emph{symmetric}
    if $\form x y=\form y x$ for all $x,y\in S$.
\end{definition}

Density on the right is equivalent to requiring the sup-lattice
homomorphism $\form 1
-:L\rightarrow 2$ to be dense, which justifies our terminology.
It is
also equivalent to requiring $y=0$ whenever $z\perp y$ for all
$z\in L$,
\ie, whenever $\form z y=\form z 0$ for all $z\in L$, which
shows that
faithfulness on the right is a stronger condition than density on the
right. Of course, these two conditions would be equivalent if we were
dealing with forms on ring modules, and equivalent to saying that
a form
is non-degenerate on the right. Hence, for sup-lattices there are 
two natural notions of non-degeneracy on the right. We
shall need both of them, so we have decided to use a
different word for each, and non-degeneracy for none in order to
avoid
ambiguity. Similar remarks apply to density and faithfulness on
the left.

In view of these remarks it may seem surprising that we have
defined
non-singular to mean the same as faithful, since for ring
modules a
non-degenerate form is not necessarily non-singular, but
(\ref{cor:nonsingular}) below shows that in the case of
sup-lattices this identification is appropriate.

\begin{definition}
Let $\varphi:L\times R\rightarrow 2$ be a 2-form, $x\in L$
and $y\in R$.
The \emph{(right) orthogonal image} of $x$ is the element
$\rorth x\in
R$ defined by
\[\rorth x=\V\{y\in R\st x\perp y\}\;.\]
Similarly, the \emph{(left) orthogonal image} of $y$ is given by
\[\lorth y=\V\{x\in L\st x\perp y\}\;.\]
\end{definition}

The correspondence between Galois connections and 2-forms can be
summarized as follows:

\begin{proposition}\label{prop:formsvsGalois}
For any 2-form between sup-lattices $L$ and $R$,
the orthogonal images
$\rorth{(-)}:L\rightarrow R$ and $\lorth{(-)}:R\rightarrow L$
form a
Galois connection, and any Galois connection between $L$
and $R$ is
uniquely determined in this way by the 2-form whose
orthogonality
relation is given by
\[x\perp y \iff x\le \lorth y\ (\iff y\le \rorth x)\;.\]
Furthermore, we have:

\begin{enumerate}

\item The following are equivalent:

\begin{enumerate}

\item $\varphi$ is dense on the right;

\item $\rorth 1=0$;

\item $0$ is the unique element
$y\in R$ such that $\lorth y=1$;

\end{enumerate}

\item The following are equivalent:

\begin{enumerate}

\item $\varphi$ is dense on the left;

\item $\lorth 1=0$;

\item $0$ is the unique element
$x\in L$ such that $\rorth x=1$;

\end{enumerate}

\item The following are equivalent:

\begin{enumerate}

\item $\varphi$ is faithful on the right;

\item $\rorth{(-)}$ is surjective;

\item $\lorth{(-)}$ is injective;

\item $\rorth{(\lorth y)}=y$ for all $y\in R$.

\end{enumerate}

\item The following are equivalent:

\begin{enumerate}

\item $\varphi$ is faithful on the left;

\item $\rorth{(-)}$ is injective;

\item $\lorth{(-)}$ is surjective;

\item $\lorth{(\rorth x)}=x$ for all $x\in L$.

\end{enumerate}
\end{enumerate}
\end{proposition}

\begin{corollary}\label{cor:nonsingular}
A 2-form between sup-lattices $L$ and $R$ is faithful if and only if
$\rorth{(-)}$ (equiv., $\lorth{(-)}$) is an antitone order
isomorphism.
\end{corollary}

Let us see some explicit examples of Galois connections in the
language
of 2-forms.

\begin{example}\label{exm:tspace}
    Let $X$ be a topological space, with topology $\opens X$.
    Then we define a 2-form between
    $\pwset X$ and $\opens X$
    by
    \[S\perp U\iff S\cap U=\emptyset\;,\]
    which is faithful on the right and dense on the left.
    More generally, given any closure operator
    $j:L\rightarrow L$ on a sup-lattice $L$, the assignment
    $x\mapsto j(x)$ defines a surjective sup-lattice homomorphism
    $j:L\rightarrow L_{j}$ (\cf\ \S\ref{sec:background}),
    and we may define a 2-form on $L\times\opp{L_{j}}$ by
    $x\perp y\iff
    x\leq y$ in $L$. Furthermore, this form is necessarily faithful
    on the right,
    and it is dense on the left
    if and only if $j(0)=0$ (\ie, the closure is dense).
    In particular,
    if we take $L$
    to be $\pwset X$ and $L_{j}$ to be the lattice of
    closed sets of
    $X$ then $\opp{L_{j}}\cong\opens X$ and we obtain
    the same as before.
\end{example}

\begin{example}
    Let $\rho$ be a binary relation between two sets $S$ and $T$.
    Then we
    have a 2-form between $\pwset S$ and $\pwset T$ given by
    \[X\perp Y\iff x\rho y\textrm{ for all }x\in X\textrm{, }y\in
    Y\;.\]
    For instance, if $S=T$ and we take $x\rho y$ to be $x\neq y$,
    we obtain
    \[X\perp Y\iff X\cap Y=\emptyset\;,\]
    as in the topological example above.
\end{example}

\begin{example}
    Let $R$ be a commutative ring, $M$ and $N$
    $R$-modules, and $f:M\times N\rightarrow R$ a bilinear 
    form. Then $f$ induces an 
    orthogonality relation between $M$ and $N$, with respect to which 
    we can define a sup-lattice 2-form
    $\varphi:\sub_{R}(M)\times\sub_{R}(N)\rightarrow 2$
    as in the previous example: for all submodules $X\subseteq M$ and
    $Y\subseteq N$,
    put $X\perp Y$ if and only if
    $f(x,y)=0$ for all $x\in 
    X$ and $y\in Y$.
    This sup-lattice 2-form is dense on the left
    (resp.\ on the right)
    if and only if the bilinear 
    form $f$ is non-degenerate on the left (resp.\ on the right).
\end{example}

\begin{definition}
    Let $\varphi:L\times R\rightarrow 2$ and $\varphi':L'\times
    R'\rightarrow 2$ be 2-forms. An \emph{orthomorphism}
    $\varphi\rightarrow\varphi'$ is a pair of sup-lattice
    homomorphisms $f:L\rightarrow L'$ and $g:R\rightarrow R'$
    such that
    for all $x\in L$ and $y\in R$ we have
    \[\form {f(x)}{g(y)} = \form x y\;.\]
    If both $f$ and $g$ are surjective the orthomorphism $(f,g)$
    is said to 
    be a \emph{quotient 
    orthomorphism}. In that case $\varphi'$ is an
    \emph{orthoquotient}, or simply \emph{quotient}, of $\varphi$.
\end{definition}

\begin{proposition}\label{prop:pseudo}
    Let $\varphi:L\times R\rightarrow 2$ and
    $\varphi':L'\times R'\rightarrow
    2$ be 2-forms, and $(f,g):\varphi\rightarrow\varphi'$ an
    orthomorphism. If $g$ is surjective then
    $(f,g)$ commutes with
    $\rorth{(-)}$ in the sense that
    \[g(\rorth x)=\rorth{f(x)}\]
    for all $x\in L$. If furthermore $f$ is strong then $(f,g)$
    preserves density on the right (\ie, $\varphi'$
    is dense on the right if $\varphi$ is), and, if
    $g$ is also dense, then
    $\varphi$ is dense on the right if
    and only if $\varphi'$ is dense on the right.
    Obvious dual statements apply to $\lorth{(-)}$
    and density on the left if $f$ and $g$
    are interchanged.
\end{proposition}

\begin{proof}
    Assume that $g$ is surjective, and let $x\in L$. Then,
    \[g(\rorth x)=g(\V\{y\st x\perp y\})=\V\{g(y)\st x\perp y\}=\]
    \[=\V\{g(y)\st f(x)\perp g(y)\}=\V\{z\in R'\st f(x)\perp
    z\}=\rorth{f(x)}\;.\]
    If furthermore $f$ is strong and $\varphi$ is dense on
    the right we have
    \[\rorth{1_{L'}}=\rorth{f(1_{L})}=g(\rorth{1_{L}})=g(0)=0\;,\]
    \ie, $\varphi'$ is dense on the right; and if
    $g$ is also dense we have $\rorth{1_{L}}=0$ if and
    only if $g(\rorth{1_{L}})=0$, and thus
    $\varphi$ is dense on the right if and only if $\varphi'$ is.
    The dual facts, with $f$ and $g$ interchanged, are proved
    in a similar way. \qed
\end{proof}

\begin{proposition}\label{prop:surjvsembed}
    Let $\varphi:L\times R\rightarrow 2$ and
    $\varphi':L'\times R'\rightarrow 2$ be 2-forms, and
    $(f,g):\varphi\rightarrow \varphi'$ an orthomorphism.
    If $g$ is surjective and $\varphi$ is faithful on the left
    then $f$ is an order embedding,
    and if $f$ is surjective and $\varphi$ is faithful
    on the right then $g$ is an order embedding.
\end{proposition}

\begin{proof}
    Assume that $g$ is surjective and
    $\varphi$ is faithful on the left, \ie,
    $\lorth{(\rorth z)}=z$ for all $z\in L$
    [\cf\ (\ref{prop:formsvsGalois})].
    Then $(f,g)$
    preserves $\rorth{(-)}$,
    by the previous proposition, and thus
    we have, for all $x,z\in L$,
    \[f(x)\le f(z)\Rightarrow f(x)\le
    \lorth{(\rorth{f(z)})}\iff f(x)\le
    \lorth{(g(\rorth z))}\iff
    f(x)\perp g(\rorth z)\]
    \[\iff x\perp \rorth z\iff x\le \lorth{(\rorth{z})}=z\;,\]
    \ie, $f$ is an order embedding.
    For $f$ surjective and $\varphi$ faithful on the right
    everything is similar.
    \qed
\end{proof}
    
\begin{corollary}
    Let $\varphi$ and $\varphi'$ be 2-forms, and
    $(f,g):\varphi\rightarrow \varphi'$ a quotient orthomorphism.
    If $\varphi$ is faithful then both $f$ and $g$ are
    order isomorphisms.
\end{corollary}

Hence, the faithful 2-forms are ``simple'' in 
the sense that their only quotient orthomorphisms are isomorphisms.
The next proposition 
states that from any 2-form we can always obtain a faithful one by 
means of a quotient, and corresponds to the fact
that any Galois connection restricts to 
a dual isomorphism between the lattices of ``closed elements''.

\begin{proposition}\label{prop:simpleorthoquotient}
    Let $\varphi:L\times R\rightarrow 2$ be a 2-form, and let
    $f:L\rightarrow L$ and $g:R\rightarrow R$ be the closure
    operators
    defined by
    $x\mapsto\lorth{(\rorth x)}$ and $y\mapsto\rorth{(\lorth y)}$.
    Let $L'=\lorth R=\{\lorth{(\rorth x)}\st x\in L\}$ and
    $R'=\rorth L=\{\rorth{(\lorth y)}\st y\in R\}$ be the
    corresponding quotients 
    of $L$ and $R$, and define a map $\varphi':L'\times 
    R'\rightarrow 2$ to be the restriction of $\varphi$ to $L'\times 
    R'$. Then 
    $\varphi'$ is a faithful 2-form and the pair $(f,g)$ defines a
    (quotient) orthomorphism from $\varphi$ to $\varphi'$.
\end{proposition}

\begin{proof}
    First we remark that for all $x\in L$ and $y\in R$ we have
    \[x\perp y\iff y\le \rorth x=\rorth{(\lorth{(\rorth x)})}
    \iff\lorth{(\rorth x)}\perp y\;.\]
    Hence, for each subset $X\subseteq L'$ and each $y\in R$ we have
    \[\lorth{(\rorth{(\V X)})}\perp y\iff\V X\perp y\iff x\perp 
    y\ \textrm{ for all }x\in X\;,\]
    which shows that $\varphi'$ preserves joins in the left variable, 
    because the join of $X$ in $L'$ is $\lorth{(\rorth{(\V X)})}$. 
    Similarly, $\varphi'$ preserves joins on the right and is thus a
    2-form, obviously 
    faithful, see (\ref{prop:formsvsGalois}).
    Finally, we also obtain, for all $x\in L$ and 
    $y\in R$,
    \[x\perp y\iff \lorth{(\rorth x)}\perp y\iff\lorth{(\rorth x)}
    \perp 
    \rorth{(\lorth y)}\;,\]
    which means $(f,g)$ is an orthomorphism. \qed
\end{proof}

\begin{definition}
    We refer to the faithful 
    2-form $\varphi'$ of the previous proposition as the
    \emph{orthogonal quotient} of $\varphi$.
\end{definition}
    
We conclude this section with the following proposition, which
will not be needed elsewhere in this paper, but which can be regarded
as the ``soft" version of (\ref{prop:surjvsembed}):

\begin{proposition}
    Let $\varphi:L\times R\rightarrow 2$ and $\varphi':L'\times R'
    \rightarrow
    2$ be 2-forms, and $(f,g):\varphi\rightarrow\varphi'$ an
    orthomorphism. If $g$ is strong and $\varphi$ is dense on the left
    then $f$ is
    dense. If $f$ is strong and $\varphi$ is dense on the right
    then $g$ is dense.
\end{proposition}

\begin{proof}
    Assume that $\varphi$ is dense on the left. If $g$ is strong
    then we have
    \[f(x)=0\Rightarrow f(x)\perp 1\iff f(x)\perp g(1)\iff x
    \perp 1\iff x=0\;,\]
    \ie, $f$ is dense. The second part of the proof is similar.
    \qed
\end{proof}

\section{Quantales and 2-forms}\label{sec:quantalesandtwoforms}

The multiplication of a quantale $Q$ has the property that,
    for all
$X\subseteq Q$ and $a\in Q$, $(\V X)\then a=0$ if and only
    if $x\then
a=0$ for all $x\in X$, and $a\then\V X=0$ if and only if
    $a\then x=0$ for
all $x\in X$. Hence, we obtain a 2-form from any quantale,
    as follows:

\begin{definition}
For any quantale $Q$, we define a 2-form $\Phi(Q)$ between $\ls(Q)$
and
$\rs(Q)$ by putting, for each $a\in\ls(Q)$ and $b\in\rs(Q)$,
\[a\perp b\iff a\then b=0\;.\]
\end{definition}

Now we study a converse to this, \ie, a way of obtaining a quantale
from
a 2-form, after which we relate the two constructions.

\begin{definition}
    Let $\varphi:L\times R\rightarrow 2$ and
    $\varphi':L'\times R'\rightarrow 2$ be
    2-forms.
    A \emph{continuous map} from $\varphi$ to $\varphi'$ is a pair
    $(f,g)$ of
    contravariant
    sup-lattice homomorphisms, where $f:L\rightarrow L'$
    and $g:R'\rightarrow R$, such that the following
    \emph{continuity} condition is satisfied: for all $x\in L$
    and
    $y\in R'$,
    \[\form{f(x)} y =\form x{g(y)}\;.\]
\end{definition}

\begin{example}\label{exm:continuousmaps}
    The above terminology is justified as follows.
    Let $X$ and $Y$ be topological spaces, with topologies $\opens X$ 
    and $\opens Y$, let $f: X\rightarrow Y$ be a map (not 
    necessarily continuous), and let
    $g:\opens Y\rightarrow\opens X$ be a sup-lattice homomorphism. 
    Seeing $X$ and $Y$ as 2-forms as in 
    (\ref{exm:tspace}), the pair $(\tilde f,g)$, where $\tilde f$ 
    is the direct image map of $f$, $\tilde f(S)=\{f(x)\st x\in S\}$,
    is a continuous
    map of 2-forms if and only if $f:X\rightarrow Y$ is a continuous
    map of 
    topological spaces and $g=f^{-1}$.
    
    A generalization of this situation can be obtained from any
    pair of
    closure operators $j$ 
    and $j'$ on sup-lattices $L$ and $L'$, respectively.
    From these 
    we obtain 2-forms on $L\times\opp{L_{j}}$ and 
    $L'\times\opp{L'_{j'}}$, as in the second part of 
    (\ref{exm:tspace}), and if $f:L\rightarrow L'$ is a 
    sup-lattice homomorphism, then $(f,g)$ is a continuous map of 
    2-forms if and only if $f$ satisfies the condition
    $f\circ j\le j'\circ f$ (\ie, $f$ is continuous with respect to 
    the closure operators) and $g$ is the restriction to
    $L_{j'}$ of 
    the right adjoint $f_{*}$---see~(\ref{prop:uniquedet}) below.
\end{example}

\begin{proposition}\label{prop:continuities}
The continuity condition is equivalent to each of the following:
\begin{enumerate}
\item $g_*(\rorth x)=\rorth{f(x)}$ for all $x\in L$,
\item $f_*(\lorth y)=\lorth{g(y)}$ for all $y\in R$,
\item $f(x)\perp g_*(\rorth x)$ and $x\perp g(\rorth{f(x)})$
    for all $x\in L$,
\item $f_*(\lorth y)\perp g(y)$ and $f(\lorth{g(y)})\perp y$
    for all $y\in R$.
\end{enumerate}
\end{proposition}

\begin{proof}
1. Continuity can be rewritten as
\[f(x)\perp y\iff x\perp g(y)\;,\]
which in turn is equivalent to
 the condition
\[y\le\rorth{f(x)}\iff g(y)\le\rorth x\;,\]
whose right-hand side is equivalent to $y\le g_*(\rorth{x})$.
\par 2. This is similar to the previous case, once we rewrite the
continuity condition as
\[f(x)\le\lorth y\iff x\le\lorth{g(y)}\;,\]
since now the left hand side is equivalent to $x\le f_*(\lorth y)$.
\par 3. From the first condition, continuity is equivalent to the
conjunction
\[g_*(\rorth x)\le\rorth{f(x)}\mbox{ and }g_*(\rorth x)
\ge\rorth{f(x)}\;.\]
The inequality $g_*(\rorth x)\le \rorth{f(x)}$ is equivalent to
$f(x)\perp g_*(\rorth x)$, and the other inequality is equivalent to
$\rorth x\ge g(\rorth{f(x)})$, \ie, to $x\perp g(\rorth{f(x)})$.
\par 4. Similar to the previous case, now using the second condition.
\qed
\end{proof}

\begin{corollary}
    Let $\varphi:L\times R\rightarrow 2$ and $\varphi':L'\times 
    R'\rightarrow 2$ be 2-forms.
    \begin{enumerate}
	\item If $\varphi$ is faithful 
    on the right then for each sup-lattice homomorphism 
    $f:L\rightarrow L'$ there is at most one sup-lattice homomorphism
    $g:R'\rightarrow R$ such that $(f,g)$ is a continuous map of 
    2-forms $\varphi\rightarrow\varphi'$.
    \item If $\varphi'$ is faithful on 
    the left then for each sup-lattice homomorphism 
    $g:R'\rightarrow R$ there is at most one sup-lattice homomorphism
    $f:L\rightarrow L'$ such that $(f,g)$ is a continuous map of 
    2-forms $\varphi\rightarrow\varphi'$.
    \end{enumerate}
\end{corollary}

\begin{proof}
    1. $\varphi$ is faithful on the right if and only if the map 
    $\rorth{(-)}:L\rightarrow R$ is surjective. Hence, the first
    condition of the 
    proposition, $g_*(\rorth x)=\rorth{f(x)}$ for all $x\in L$, 
    completely determines the right adjoint $g_{*}$, and thus it
    determines $g$.\\
    2. Similar, taking into account the second condition of the 
    proposition. \qed
\end{proof}

Notice that we do not state, \eg\ in the first part of this corollary, 
that for every $f$ there is a $g$ such that $(f,g)$ is continuous.
The second part of (\ref{exm:continuousmaps})
provides an example in which for only some $f$ this holds, and
the following proposition gives a necessary and suficient
condition for such $g$ to exist.

\begin{proposition}\label{prop:uniquedet}
    Let $\varphi:L\times R\rightarrow 2$ and $\varphi':L'\times 
    R'\rightarrow 2$ be 2-forms, both faithful on 
    the right. Let also $f:L\rightarrow L'$ be a
    homomorphism. Then the following are equivalent:
    \begin{enumerate}
    \item There is a homomorphism $g:R'\rightarrow R$ such that
    $(f,g)$ 
    is continuous;
    \item $f(\lorth{(\rorth 
    x)})\le\lorth{(\rorth{f(x)})}$ for all $x\in L$ (\ie, $f$ is 
    continuous with respect to the closure operators
    $\lorth{(\rorth{(-)})}$.
    \end{enumerate}
    If in addition $\varphi$ is faithful on the left, then
    for each $f:L\rightarrow L'$ 
    there is exactly one $g:R'\rightarrow R$ such that $(f,g)$ is 
    continuous.
\end{proposition}

\begin{proof}
    Assume that (1) holds.
    Then we have $f(x)\perp \rorth{f(x)}$, and
    \[f(x)\perp \rorth{f(x)}\iff x\perp g(\rorth{f(x)})\iff
    \lorth{(\rorth x)}\perp g(\rorth{f(x)})\iff\]
    \[\iff f(\lorth{(\rorth x)})\perp\rorth{f(x)}\iff
    f(\lorth{(\rorth 
    x)})\le\lorth{(\rorth{f(x)})}\;.\]
    Now assume that (2) holds. Write $j$ for the closure operator
    $\lorth{(\rorth{(-)})}$ on $L$, and $k$ for the similar closure 
    on $L'$. Then (2) is the condition $f\circ j\le k\circ f$. We 
    shall prove that the image of the restriction of $f_{*}$ to 
    $L'_{k}$ is contained in $L_{j}$. Indeed, this is equivalent to 
    the condition that $j\circ f_{*}\circ k\le f_{*}\circ k$, which 
    holds because
    \[f\circ j\le k\circ f\iff j\le f_{*}\circ k\circ f\Rightarrow
    j\circ f_{*}\circ k\le f_{*}\circ k\circ f\circ f_{*}\circ k\le
    f_{*}\circ k\;,\]
    where the latter inequality follows from the fact that $f\circ 
    f_{*}\le{\rm id}_{L'}$ and $k\circ k=k$. Hence, $f_{*}$ defines
    a meet preserving 
    map $L'_{k}\rightarrow L_{j}$. Due to right faithfulness of 
    $\varphi$ and $\varphi'$ we have order isomorphisms
    $L_{j}\cong\opp R$ and $L'_{k}\cong\opp{R'}$, and thus 
    $g:R'\rightarrow R$ can be defined by composing $f_{*}$ with 
    the isomorphisms. Finally, if $\varphi$ is also faithful on the 
    left we have $\lorth{(\rorth x)}=x$ for all $x\in L$, and thus 
    $f$ is trivially continuous with respect to
    $\lorth{(\rorth{(-)})}$. \qed
\end{proof}

Clearly, continuous maps are closed under composition,
and thus we obtain
another category of 2-forms, which furthermore is sup-lattice
enriched.
In particular, then, the continuous endomaps of any 2-form
form a unital
quantale:

\begin{definition}
Let $\varphi:L\times R\rightarrow 2$ be a 2-form. The
\emph{quantale of
$\varphi$}, denoted by $\Q(\varphi)$, is the quantale
of continuous
endomaps of $\varphi$, with
$(f,g)\le (f',g')$ if and only if $f(x)\le f'(x)$ and $g(y)\le g'(y)$
for all $x\in L$ and $y\in R$,
and with multiplication given by
$(f,g)\then (f',g')=(f'\circ f, g\circ g')$.
\end{definition}

In the case of a symmetric 2-form $\varphi$ we have
$(f,g)\in\Q(\varphi)$ if and only if $(g,f)\in\Q(\varphi)$,
and at once we remark:

\begin{proposition}\label{prop:symvsinv}
Let $\varphi:L\times L\rightarrow 2$ be a 
	symmetric
	2-form. Then
	the quantale $\Q(\varphi)$ is involutive,
         with the 
	involution given by $(f,g)^{*}=(g,f)$.
         Conversely,
         if $Q$ is an involutive quantale then $\Phi(Q)$ is
    isomorphic 
	to a symmetric 2-form, and $\Q(\Phi(Q))$ is involutive,
    with the involution given by $(f,g)^{*}=(g',f')$, where
	$f':\rs(Q)\rightarrow\rs(Q)$ and $g':\ls(Q)\rightarrow\ls(Q)$
    are 
	defined by
	$f'(y)=f(y^{*})^{*}$ and $g'(x)=g(x^{*})^{*}$.
\end{proposition}

\begin{example}\label{exm:endomorphisms}
    Let us relate the quantales of endomorphisms of
    2-forms to the
    well known endomorphism quantales of
    sup-lattices.
    \begin{enumerate}
	\item Let $\varphi:L\times R\rightarrow 2$ be a
    faithful 2-form. 
	From (\ref{prop:uniquedet}) it follows that
	the quantales $\Q(\varphi)$ and $\Q(L)$ are isomorphic.
	\item\label{exm:endomorphisms.2} Let $L$ be a sup-lattice,
    and define a 2-form
	$\varphi:L\times\opp L\rightarrow 2$ by $x\perp y$
    if and only if 
	$x\le y$ in $L$ (in other words, consider the Galois
    connection 
	between $L$ and $\opp L$
	defined by the identity
    $\rorth{(-)}={\rm id}_{L}:L\rightarrow \opp{(\opp L)}$).
    This 2-form is 
	faithful, and thus the quantales
	$\Q(\varphi)$ and $\Q(L)$ are isomorphic.
	\item Let
	$\varphi:L\times L\rightarrow 2$ be both
    symmetric and faithful. 
	Then $\Q(\varphi)$ is isomorphic to $\Q(L)$,
    which is thus involutive.
	The involution is defined on $\Q(L)$ in the usual 
	way for quantales of endomorphisms on self-dual
	sup-lattices~\cite{MulvPell92}:
	\[f^{*}(y)=\rorth{(\V\{x\st f(x)\le\rorth y\})}\;.\]
\end{enumerate}
\end{example}

\begin{example}
Kruml~\cite{KrumPhD} defines a
\emph{Galois quantale} to 
be a quantale
$\Q(G)=\{(f,g)\in\Q(S)\times \Q(T)\st g\circ G=G\circ f\}$
for some sup-lattice homomorphism $G:S\rightarrow T$.
From
(\ref{prop:continuities}) it follows that Galois 
quantales
are the same as quantales of 2-forms: $\Q(G)$ is
isomorphic to
$\Q(\varphi)$ for the 2-form $\varphi:S\times
\opp T\rightarrow 2$ 
such that $\rorth{(-)}=G$.
\end{example}

\begin{lemma}\label{lem:sidedisos}
Let $\varphi:L\times R\rightarrow 2$ be a dense two-form. Then
$\ls(\Q(\varphi))$ is order isomorphic to $L$, and
$\rs(\Q(\varphi))$ is
order isomorphic to $R$. Furthermore $\Q(\varphi)$
is a factor 
quantale.
\end{lemma}

\begin{proof}
First we remark that $\Q(\varphi)$ is a subquantale of
$\Q(L)\times\Q^*(R)$, where $\Q^{*}(R)$ is the 
quantale $\Q(R)$ with reversed multiplication,
\ie, with $f\then g=f\circ 
g$. Also, the top of $\Q(L)\times\Q^*(R)$
belongs to $\Q(\varphi)$
because $\varphi$ is dense: for all $x\in L$ and $y\in R$,
if either $x=0$ or $y=0$ then both conditions
$1_{\Q(L)}(x)\perp y$ and $x\perp 1_{\Q^*(R)}(y)$ are
true, whereas if $x\neq 0$ and $y\neq 0$ then
both conditions are false.
Hence, the left-sided
elements of $\Q(\varphi)$ are precisely those which
are left-sided as elements of
$\Q(L)\times\Q^*(R)$, \ie, they are the continuous maps
of the form
$(c_l,a_r)$ for some $l\in L$ and $r\in R$, where $c_{l}$
and $a_{r}$ 
are respectively a ``constant'' map and an annihilator,
as described 
in \S\ref{sec:background}.
Hence, continuity means that for
any pair of elements $x\in L$ and $y\in R$ we must have
\[\form{c_l(x)}y=\form x{a_r(y)}\;.\]
Taking $x=1$ yields
$\form{l}y=\form 1{a_r(y)}$,
and thus
\[\begin{array}{rcll}
l\perp y &\iff& 1\perp a_r(y)\\
&\iff& a_r(y)=0 & \textrm{(due to density)}\\
&\iff& y\le r\;.
\end{array}\]
Therefore a necessary condition for continuity is
$r=\rorth l$. The condition is also sufficient, since:
\begin{itemize}
\item if $x=0$ then, trivially, $\form{c_l(x)}y=\form x{a_r(y)}=0$;
\item if $x\neq 0$ then $c_l(x)=l$ and we obtain
\[\begin{array}{rcll}
\form{c_l(x)}y=0 &\iff& l\perp y\\
&\iff&  y\le r\\
&\iff& a_r(y)=0\\
&\iff& \form x{a_r(y)}=0\;,
\end{array}\]
\end{itemize}
where the latter step follows from density on the left and the fact
that $a_r(y)$ equals either $0$ or $1$.
Hence, the
generic form of a left-sided element of $\Q(\varphi)$ is
$(c_l,a_{\rorth
l})$, which means we have an assignment
$l\mapsto(c_l,a_{\rorth l})$ that
defines a surjective map $L\rightarrow\ls(\Q(\varphi))$.
Furthermore, we
have $l\le k$ if and only if $c_l\le c_k$, and $l\le k$ implies
$a_{\rorth l}\le a_{\rorth k}$, which makes the map
$L\rightarrow\ls(\Q(\varphi))$ an order-isomorphism.
For right-sided
elements the proof is analogous: each right-sided element of
$\Q(\varphi)$ must be of the form $(a_l,c_r)$ for some
$l\in L$ and $r\in
R$, and continuity is the condition
\[\form{a_l(x)}y=\form x{c_r(y)}\;.\]
Taking again density of $\varphi$ into account, we conclude that
$l=\lorth r$, and thus the right-sided elements must be of the form
$(a_{\lorth r},c_r)$. Hence, $R$ is isomorphic to $\rs(\Q(\varphi))$. 
Finally, the only elements that are simultaneously left- and 
right-sided are those for which $(c_l,a_{\rorth l})=(a_{\lorth 
r},c_r)$,
with $l\in L$ and $r\in R$. The only solutions are
$(1,1)$ and $(0,0)$,
corresponding respectively to $l=r=1$ and $l=r=0$,
\ie, $\Q(\varphi)$ 
is a factor. \qed
\end{proof}

\begin{example}
    Let $L$ be a sup-lattice, and $\varphi$ the 2-form on 
    $L\times\opp L$ of 
    (\ref{exm:endomorphisms})-(\ref{exm:endomorphisms.2})
    [\ie, with $x\perp y\iff x\le y$]. From the isomorphism
    $\Q(L)\cong\Q(\varphi)$ we immediately
    obtain the well known 
    isomorphisms $\ls(\Q(L))\cong L$ and $\rs(\Q(L))\cong\opp 
    L$~\cite{MulvPell92}.
\end{example}

\begin{theorem}
Let $\varphi:L\times R\rightarrow 2$ be a dense 2-form. Then $\varphi$
and $\Phi(\Q(\varphi))$ are isomorphic 2-forms.
\end{theorem}

\begin{proof}
All that we have to do is show that the isomorphisms of the previous
lemma commute with the forms, \ie, that for all $l\in L$ and
$r\in R$ we
have $l\perp r$ if and only if, in $\Q(\varphi)$,
the following condition
holds,
 \[(c_l,a_{\rorth l})\then(a_{\lorth r},c_r)=(0,0)\;,\]
or, equivalently, if and only if the two following conditions hold:
(i)~$a_{\lorth r}\circ c_l=0$ and (ii)~$a_{\rorth l}\circ c_r=0$.
Since we have
$c_l(1)=l$, condition~(i) holds if and only if
$a_{\lorth r}(l)=a_{\lorth
r}(c_l(1))=0$, which is equivalent to $l\le\lorth r$.
Similarly,
condition~(ii) holds if and only if $a_{\rorth l}(r)=0$,
which is
equivalent to $r\le\rorth l$. Hence, both~(i) and~(ii)
are equivalent to
$l\perp r$. \qed
\end{proof}

It is not in general true that for a quantale $Q$ we have
$Q\cong \Q(\Phi(Q))$, but there is always a \emph{comparison 
homomorphism}
$\kappa:Q\rightarrow\Q(\Phi(Q))$:

\begin{proposition}
    Let $Q$ be a quantale, and let $a\in Q$.
    The right action of $a$ on $\ls(Q)$ and the 
    left action of $a$ on $\rs(Q)$ jointly define a
    continuous endomap
    $((-)\then a,a\then(-))$ of the 
    2-form $\Phi(Q)$. The comparison homomorphism 
    $\kappa:Q\rightarrow\Q(\Phi(Q))$ defined by
    $a\mapsto((-)\then 
    a,a\then(-))$ is a homomorphism of quantales,
    unital if $Q$ is 
    unital. If furthermore $Q$ is involutive
    [in which case 
    $\Q(\Phi(Q))$ is involutive, see
    (\ref{prop:symvsinv})]
    then $\kappa$ preserves the involution.
\end{proposition}

\begin{proof}
    The first part is immediate from the associativity of
    multiplication
    in $Q$, for:
    \[\form{x\then a}{y}=0\iff (x\then a)\then y=0\iff x\then 
    (a\then y)=0\iff\form x{a\then y}=0\;.\]
    Now let us see that $\kappa$ is a 
    homomorphism of quantales. First,
    it preserves multiplication
    because composition of continuous
    maps of 2-forms is defined by
    $(f',g')\circ (f,g)=(f'\circ f, g\circ g')$, and thus
    for all $a,b\in Q$ the product
    $\kappa(a)\then\kappa(b)$ equals
     \[\kappa(b)\circ \kappa(a)=
    (((-)\then a)\then b,a\then(b\then(-)))
    =((-)\then (a\then b), (a\then b)\then (-))=
    \kappa(a\then b)\;.\]
    If $Q$ is 
    unital then $\kappa(e)$ is the unit of $\Q(\Phi(Q))$, 
    and if $Q$ is involutive we have
    $\kappa(a^{*})(x,y)=(x\then a^{*},a^{*}\then y)=
    ((a\then x^{*})^{*},(y^{*}\then a)^{*})=
    \kappa(a)^{*}(x,y)$, see (\ref{prop:symvsinv}). \qed
\end{proof}

The comparison homomorphism is injective if and only if,
for all $a,b\in Q$, if $x\then a=x\then b$ and 
$a\then y=b\then y$ for all $x\in\ls(Q)$
and $y\in\rs(Q)$ then we have $a=b$. A quantale satisfying this 
condition is usually said to be
\emph{faithful}~\cite{Krum02,Pase97,PellRosi97}.

\section{Principal quantale modules}\label{sec:principalqumods}

\begin{definition}
    Let $Q$ be a quantale. A left $Q$-module $M$ is \emph{principal}
    if it has a \emph{generator}, \ie, an element
    $x\in M$ such that $Qx=\{ax\st a\in Q\}=M$.
    Similar definitions
    apply to right modules.
\end{definition}

Some basic obvious properties of principal modules are the following:

\begin{proposition}\label{prop:principalmodules}
    Let $Q$ be a quantale.
    \begin{enumerate}
	\item Any left $Q$-module quotient of a principal
    left $Q$-module is 
	principal.
	\item If $M$ is a principal left $Q$-module
	then it is a left $Q$-module quotient of $Q$.
	\item If $Q$ is unital then $M$ is a principal $Q$-module 
    if and only
	if it is a left $Q$-module quotient of $Q$.
    \end{enumerate}
\end{proposition}

\begin{proof}
    1. If $f:M\rightarrow N$ is a surjective homomorphism of left 
    $Q$-modules and $M$ has a generator $x$ then $f(x)$ is a 
     generator of $N$.\\
    2. If $M$ is a left $Q$-module with a  generator $x$ then 
    the map $Q\rightarrow M$ defined by $a\mapsto ax$ is a
    surjective 
    homomorphism of left $Q$-modules.\\
    3. If $Q$ is unital then $e$ is a  generator of
    itself as a 
    module. The rest follows from the previous two. \qed
\end{proof}

\begin{definition}
    Let $Q$ be a quantale, and $M$ a left $Q$-module.
    An element $x\in 
    M$ is \emph{invariant} if $ax\le x$ for all $a\in Q$
    (equivalently, if 
    $1_{Q}x\le x$).
\end{definition}

Hence, the left-sided elements of a quantale $Q$ are the
invariant elements of
$Q$ when $Q$ is seen as a left module over itself.

\begin{proposition}\label{prop:updownsegments}
    Let $Q$ be a quantale, $M$ a left $Q$-module, and $m$
    an element 
    of $M$. The following are equivalent:
    \begin{enumerate}
    \item
    $\upsegment m$ is a left $Q$-module quotient of $M$,
    with action defined by $(a,x)\mapsto ax\vee m$ and quotient
    projection
    $Q\rightarrow\upsegment m$ given by $x\mapsto x\vee m$.
    \item
    $\downsegment m$ is a left $Q$-submodule of $M$.
    \item $m$ is an invariant element of $M$.
    \end{enumerate}
\end{proposition}

\begin{proof}
    ($1\Leftrightarrow 3$) Condition 1 holds if and only if the map
    $(-)\vee m:M\rightarrow M$ is a nucleus of left $Q$-modules.
    So assume that $m$ is invariant. Then, for all $a\in Q$ and $x\in 
    M$ we have
    \[a(x\vee m)=ax\vee am\le ax\vee m\;,\]
    \ie, $(-)\vee m$ is a nucleus.
    Now assume that $(-)\vee m$ is a nucleus.
    Then for all $a\in Q$ we 
    have $am=a(0\vee m)\le a0\vee m=m$, \ie, $m$ is invariant.\\
    ($2\Leftrightarrow 3$) $\downsegment m$ is a sub-sup-lattice,
    so it is a submodule if and 
    only if it is closed for the action. Let then $m$ be invariant, 
    $x\in\downsegment m$, and $a\in Q$. Then $ax\le am\le m$. Let 
    now $\downsegment m$ be a submodule. Then $am\in \downsegment m$ 
    for all $a\in Q$, \ie, $m$ is invariant. \qed
\end{proof}

\begin{example}\label{exm:subM}
    Let $R$ be a ring, and $M$ a left $R$-module. Then $\sub(M)$ is a 
    left module over the quantale $\sub(R)$, and an invariant element 
    $N\in\sub(M)$ is the same as a submodule of $M$. Given such a 
    submodule $N$, it follows that $\sub(N)$ is a left
    $\sub(R)$-submodule of 
    $\sub(M)$ and coincides with $\downsegment N$, whereas 
    $\upsegment N$, which is order-isomorphic to $\sub(M/N)$,
    is also 
    isomorphic as a left $\sub(R)$-module when $\upsegment N$ is 
    given the action of (\ref{prop:updownsegments}) and 
    $\sub(M/N)$ is given the action induced by the left 
    $R$-module structure of $M/N$.
\end{example}

\begin{proposition}\label{prop:invvsleftsided}
    Let $Q$ be a quantale, $M$ a left $Q$-module, and $x\in M$.
    The map $(-)x:Q\rightarrow M$ sends left-sided elements to 
    invariant elements, and the residuation $(-)/x:M\rightarrow Q$
    sends
    invariant elements to left-sided 
    elements. Furthermore, if $x$ is a  generator 
    then the residuation also reflects left-sided elements back into 
    invariant ones, \ie, $m\in M$ is invariant
    if and only if $m/x$ is 
    left-sided.
\end{proposition}

\begin{proof}
    If $a\in Q$ is left-sided then
    $1(ax)=(1a)x\le ax$, \ie, $ax$ is invariant.
    Now assume $m$ is invariant.
    We always have $(m/x)x\le m$, and thus
    $1(m/x)x\le 1m$. Since $m$ is invariant we also have
    $1(m/x)x\le m$, which is equivalent to $1(m/x)\le m/x$, \ie, 
    $m/x$ is left-sided. Now assume that $x$ is a  generator. 
    Then the map $(-)x$ is surjective,
    the inequality $(m/x)x\le m$ becomes the equality $m=(m/x)x$,
    and thus $m$ is the 
    invariant element to which $(-)x$ maps any left-sided element
    of the form $m/x$. \qed
\end{proof}

\begin{corollary}
    Let $Q$ be a quantale, $M$ a left $Q$-module, and $x\in M$. 
    Then $\ann(x)$ is left-sided.
\end{corollary}

\begin{proof}
    $\ann(x)=0/x$ is left-sided because $0$ is invariant. \qed
\end{proof}

\begin{proposition}\label{prop:densequotient}
    Let $Q$ be a quantale, $M$ a left $Q$-module, and $x\in M$ a 
     generator. Then the quotient 
    $(-)x:Q\rightarrow M$ factors 
    through the quotient $(-)\vee\ann(x):Q\mapsto \upsegment\ann(x)$ 
    and a dense homomorphism $\varphi:\upsegment\ann(x)\rightarrow M$.
    Furthermore, $\varphi$ restricts to an order
    isomorphism $M/x=\{m/x\st m\in M\}\cong M$.
\end{proposition}

\begin{proof}
    For each $a\in Q$ we have $ax=ax\vee 
    0=ax\vee\ann(x)x=(a\vee\ann(x))x$. Hence, 
    we have $(a\vee\ann(x))x\le ax$, which is equivalent to
    $a\vee\ann(x)\le (ax)/x$, \ie,
    the closure operator
    $a\mapsto (ax)/x$ on $Q$ is greater or equal to $a\mapsto 
    a\vee\ann(x)$, and thus $\varphi$ is just the restriction of
    $(-)x$ to $\upsegment\ann(x)$. It is dense because
    $ax=0$ is equivalent to $a\le\ann(x)$, and $M$ is isomorphic to
    $M/x$ because $\varphi$ is surjective.
    \qed
\end{proof}

\begin{corollary}
    $x$ is a  generator of $M$ if and only if 
    $(\upsegment\ann(x))x=M$.
\end{corollary}

In~\cite{MulvPell01} a notion of \emph{point} of a quantale is 
based on having ``enough  generators'', whereas in other papers 
it is related only to a kind of
irreducibility~\cite{Krum02,Pase97,PellRosi97}. Since these notions 
have some importance, we devote the rest of this section
to some simple
results relating the two, although they are not needed
in the rest 
of the paper.

\begin{definition}
    Let $Q$ be a quantale, and $M$ a left $Q$-module. $M$ is said to 
    be \emph{irreducible} if it has no invariant elements besides 
    $0$ and $1$, and \emph{everywhere principal} if for all non-zero
    $m\in M$ there is a generator $x\le m$.
\end{definition}

\begin{theorem}
    Let $Q$ be a quantale, and $M$ a left $Q$-module. If either of 
    the following two conditions holds, then $M$ is irreducible:
    \begin{enumerate}
	\item $M$ is everywhere principal.
	\item $M$ has a  generator $x$ such that $\ann(x)$ is a 
	maximal left-sided element of $Q$.
    \end{enumerate}
\end{theorem}

\begin{proof}
    1. For this it suffices to see that if $x\in M$ is any 
    generator then $1_{M}$ is the only invariant above $x$. So
    assume that $m$ is an invariant such that $x\le m$,
    where $x$ is a 
    generator.
    Then $1_{M}=\V M=\V Qx =1_{Q}x\le 1_{Q}m\le m$.\\
    2. Let $m\in M$ be 
    invariant. Then $m/x$ is left-sided, by
    (\ref{prop:invvsleftsided}),
    and thus either $m/x=\ann(x)$ 
    or $m/x=1_{Q}$. But $m=(m/x)x$ because $x$ is a  generator, 
    and thus $m=\ann(x)x=0_{M}$ or $m=1_{Q}x=\V Qx=\V M=1_{M}$, \ie, 
    $M$ is irreducible. \qed
\end{proof}

\begin{remark}\label{rem:MPpoints}
    In~\cite{MulvPell01} the \emph{points} of a(n involutive)
    quantale $Q$ are certain 
    right $Q$-modules 
    which are atomic as sup-lattices and whose atoms are  
    generators, being thus everywhere principal (see also 
    \S\ref{sec:invmods}). In certain places
    in~\cite{MulvPell01} the 
    hypothesis that the module satisfies an additional
    condition known as
    \emph{non-triviality} is assumed. Although formulated
    differently, this is equivalent
    to the requirement that $\ann(x)$ be a maximal 
    right-sided element for some atom $x$, which itself implies 
    irreducibility, by the previous theorem.
\end{remark}

\section{Quantale modules and 2-forms}
\label{sec:modsand2forms}

Let $Q$ be a quantale, $\varphi:L\times R\rightarrow 2$ a 2-form, and
$h:Q\rightarrow\Q(\varphi)$ a homomorphism of quantales.
For each $a\in 
Q$, the continuous endomap
$h(a)$ is a contravariant pair of maps that defines a right action 
of $Q$ on $L$ and a left action of $Q$ on $R$:
\[h(a)=((-)a,a(-))\;.\]
Examples of such homomorphisms are the comparison homomorphisms 
defined in \S\ref{sec:quantalesandtwoforms}, for which we have
$\varphi=\Phi(Q)$.
By definition of
continuity the actions of $Q$ on these modules 
satisfy, for all $x\in L$, $y\in R$, and $a\in Q$, the following 
``middle-linearity'' condition:
\[\form{xa}y = \form x{ay}\;.\]
(In other 
words, the 2-form can be identified with a sup-lattice homomorphism 
from
$L\otimes_{Q} R$, rather than just $L\otimes R$, to $2$---for tensor 
products of sup-lattices see~\cite{JoyaTier84}.)
In this section we shall study such pairs of modules:

\begin{definition}
    Let $Q$ be a quantale, and $\varphi:L\times R\rightarrow 2$
    a 2-form.
    An \emph{action} of $Q$ on $\varphi$ consists
    of a right action of $Q$ on $L$
    and a left action of $Q$ on $R$,\footnote{Our notation was 
    suggested by the fact that $L$ is the ``left part'' of the 
    2-form, and $R$ is the ``right part''. While unfortunately this
    has led to $L$ 
    being a \emph{right} module and $R$ a \emph{left} module, 
    the notation
    is consistent with the fact that often $L$ is
    $\ls(Q)$ and $R$ 
    is $\rs(Q)$ 
    for some quantale $Q$.} such 
    that for all $x\in L$, $y\in R$, and $a\in Q$ we have
    $\form{xa}y=\form x{ay}$. When the latter 
    condition holds we say the 2-form is 
    \emph{balanced} (with respect to the $Q$-modules $L$ and
    $R$), or that it 
    is a \emph{2-form over} $Q$.
    If $Q$ is unital, the 2-form is \emph{unital} if
    both $L$ and $R$ are unital modules.
\end{definition}

\begin{proposition}\label{prop:equivalentformsoverQ}
    Let $Q$ be a quantale, $L$ a right $Q$-module, $R$ a left 
    $Q$-module, and $\varphi:L\times R\rightarrow 
    2$ a sup-lattice 2-form. The following are equivalent:
\begin{enumerate}
\item $a\backslash(\rorth x)=\rorth{(xa)}$ for all $x\in L$ and $a\in Q$,
\item $(\lorth y)/a=\lorth{(ay)}$ for all $y\in R$
and $a\in Q$,
\item $xa\perp a\backslash(\rorth x)$ and $x\perp a(\rorth{(xa)})$ for
all $x\in 
L$ and $a\in Q$,
\item $(\lorth y)/a\perp ay$ and $(\lorth{(ay)})a\perp y$
for all $y\in R$ and $a\in Q$,
\item $\varphi$ is a 2-form over $Q$.
\end{enumerate}
\end{proposition}

\begin{proof}
    This is an immediate consequence of (\ref{prop:continuities}),
for being 
    a 2-form over $Q$ is the same as the map $((-)a,a(-))$ being 
    continuous, and $(-)/a$ and $a\backslash(-)$ are the right 
    adjoints to $(-)a$ and $a(-)$, respectively. \qed
\end{proof}

Now we introduce the notion of orthomorphism that is appropriate in 
the present context.

\begin{definition}
    Let $\varphi$ and $\varphi'$ be 2-forms over a quantale $Q$. An 
    \emph{orthomorphism over $Q$},
    or simply a \emph{$Q$-orthomorphism}, is an
    orthomorphism $(f,g):\varphi\rightarrow\varphi'$ such that $f$
    is a homomorphism of right 
    $Q$-modules and $g$ is a homomorphism of left $Q$-modules.
\end{definition}

Balance is 
preserved by surjections, as follows:

\begin{lemma}
    Let $Q$ be a quantale, and $\varphi:L\times R\rightarrow 2$ a 
    2-form over $Q$. Let also $\varphi':L'\times R'\rightarrow 2$ be 
    any 2-form, such that $L'$ is a right $Q$-module, and $R'$ is a 
    left $Q$-module. Let $(f,g):\varphi\rightarrow\varphi'$ be an
    orthomorphism 
    such that both $f$ and $g$ are surjective $Q$-module 
    homomorphisms (resp.\ right and left). Then $\varphi'$ is 
    balanced.
\end{lemma}

\begin{proof}
    Let $x'\in L'$, $y'\in R'$, and $a\in Q$. Due to surjectivity, 
    there is $x\in L$ such that $x'=f(x)$, and $y\in R$ such that
    $y'=g(y)$. Hence,
    \[\varphi'(x'a,y')=\varphi'(f(x)a,g(y))=\varphi'(f(xa),g(y))=\]
    \[=\varphi(xa,y)=\varphi(x,ay)=\cdots=\varphi'(x',ay')\;.\qed\]
\end{proof}

\begin{lemma}
    Let $Q$ be a quantale, and
    $\varphi:L\times R\rightarrow 2$ a 2-form over $Q$. Then the 
    closure operators on $L$ and $R$ defined by
    $x\mapsto\lorth{(\rorth x)}$ and $y\mapsto\rorth{(\lorth y)}$
    are nuclei of $Q$-modules.
\end{lemma}

\begin{proof}
We prove this only for $\lorth{(\rorth{(-)})}$, as the other case is 
similar. Let $x\in L$ and $a\in Q$. The 
condition
$xa\perp \rorth{(xa)}$ is always true and equivalent to
$x\perp a(\rorth{(xa)})$, which in turn is equivalent to
$a(\rorth{(xa)})\le\rorth x=\rorth{(\lorth{(\rorth x)})}$, \ie, to
$\lorth{(\rorth x)}\perp a(\rorth{(xa)})$. Finally, this is 
equivalent to $(\lorth{(\rorth x)})a\perp \rorth{(xa)}$, \ie,
$(\lorth{(\rorth x)})a\le\lorth{(\rorth{(xa)})}$, which is 
precisely the statement that $\rorth{(\lorth{(-)})}$ is a nucleus 
of right $Q$-modules. \qed
\end{proof}

\begin{theorem}\label{thm:simpleQquotient}
    Let $Q$ be a quantale, $\varphi:L\times R\rightarrow 2$ a 2-form 
    over $Q$,
    and
    $f:L\rightarrow L$ and $g:R\rightarrow R$ the closure operators
    defined by
    $x\mapsto\lorth{(\rorth x)}$ and $y\mapsto\rorth{(\lorth y)}$.
    Let $L'=\lorth R=\{\lorth{(\rorth x)}\st x\in L\}$ and
    $R'=\rorth L=\{\rorth{(\lorth y)}\st y\in R\}$ be the
    corresponding quotients 
    of $L$ and $R$, and let the 2-form $\varphi':L'\times 
    R'\rightarrow 2$ be the restriction of $\varphi$ to $L'\times 
    R'$. Then $\varphi'$ is a 2-form over $Q$.
\end{theorem}
[In other words, if $\varphi$ is a 2-form over $Q$ then its simple 
quotient, as defined in (\ref{prop:simpleorthoquotient}),
is also a 2-form over $Q$.]

\begin{proof}
    Corollary of the previous lemmas and
    (\ref{prop:simpleorthoquotient}). \qed
\end{proof}

\begin{lemma}
    Let $Q$ be a quantale, and $n\in Q$.
    Define a map $\varphi_{n}:Q\times Q
    \rightarrow 2$ by
    \[\varphi_{n}(x,y)=0\iff x\then y\le n\;.\]
    Then $\varphi_{n}$ is a 2-form, and it is balanced with respect
    to $Q$, seen 
    both as a right and a left module over itself.
\end{lemma}

\begin{proof}
    We have $\varphi_{n}(\V_{i} x_{i}, y)=0$ if and only if
    $\V_{i} x_{i}\then y\le n$, which holds if and only if
    $x_{i}\then y\le n$ for all $i$, \ie,
    $\varphi_{n}(x_{i}, y)=0$ for 
    all $i$. A similar fact holds for joins on the right variable, 
    and thus $\varphi_{n}$ is a 2-form on $Q\times Q$.
    Furthermore, this 
    2-form is obviously balanced with respect to the
    actions of $Q$ on itself, 
    due to the associativity of the multiplication in $Q$. \qed
\end{proof}

\begin{definition}
    Let $Q$ be a quantale. A 2-form $\varphi:L\times R\rightarrow 2$
    over $Q$ is \emph{principal} if both 
    $L$ and $R$ are principal $Q$-modules.
\end{definition}

Hence, if $Q$ is unital and $n\in Q$ then $\varphi_{n}$ is 
principal, and so is any of its quotients.

\begin{definition}
    Let $Q$ be a quantale, and $\varphi:L\times R\rightarrow 2$
    a principal 
    2-form with  generators $x\in L$ and $y\in R$. The
    \emph{orthogonalizer} of $x$ and $y$ is defined to be
    \[\orth(x,y)=\V\{a\in Q\st x\perp ay\}\;.\]
\end{definition}
Notice that
the following equivalences hold,
\[a\le\orth(x,y)\iff x\perp ay\iff ay\le\rorth x\iff 
a\le(\rorth x)/y\;,\]
and thus $\orth(x,y)=(\rorth x)/y$. Also, we have
\[x\perp ay\iff xa\perp y\iff xa\le\lorth y\;,\]
whence $\orth(x,y)=x\backslash(\lorth 
y)$.

\begin{theorem}\label{thm:principalforms}
    Let $Q$ be a quantale, $\varphi:L\times R\rightarrow 2$ a 
    principal 2-form over $Q$, and $n=\orth(x,y)$ for a pair of  
    generators $x\in L$ and $y\in R$. Then 
    there is a $Q$-orthoquotient
    $(f,g):\varphi_{n}\rightarrow\varphi$.
\end{theorem}

\begin{proof}
    The 
    map $f:Q\rightarrow L$ defined by
    $a\mapsto xa$ is a surjective right $Q$-module homomorphism,
    and the map
    $g:Q\rightarrow R$ defined by $b\mapsto by$ is a
    surjective left 
    $Q$-module homomorphism. Finally, for all $a,b\in Q$
    we have
    \[\varphi_{n}(a,b)=0\iff a\then b\le n\iff
    x\perp aby\iff xa\perp 
    by\iff f(a)\perp g(b)\;,\]
    which shows that $(f,g)$ is a $Q$-orthomorphism. \qed
\end{proof}

In particular, this gives us a classification of the 
principal 2-forms over a unital quantale $Q$:

\begin{corollary}
    Let $Q$ be a unital quantale. Then $\varphi$ is
    a principal 2-form over $Q$ if and only if
    it is a $Q$-orthoquotient of $\varphi_{n}$ for some $n\in Q$.
\end{corollary}

If furthermore the 2-forms are required to be faithful we obtain:

\begin{corollary}
    Let $Q$ be a unital quantale. Then $\varphi$ is
    a faithful principal 2-form over $Q$ if and only if
    it is isomorphic to the orthogonal quotient
    [necessarily a $Q$-orthoquotient, by
    (\ref{thm:simpleQquotient})]
    of $\varphi_{n}$ for some $n\in Q$.
\end{corollary}

We conclude this section relating these 2-forms with the upsegment 
modules of \S\ref{sec:principalqumods}.

\begin{lemma}\label{lemma:denseforms}
    Let $Q$ be a quantale, $n\in Q$, $r\in\rs(Q)$, and $l\in\ls(Q)$, 
    such that $r\vee l\le n$, and let $\psi:\upsegment 
    r\times\upsegment l\rightarrow 2$ be the
    restriction of $\varphi_{n}$ to $\upsegment r\times\upsegment 
    l$. Then:
    \begin{enumerate}
    \item with the quotient module structures of $\upsegment r$ and 
    $\upsegment l$, $\psi$ is a 2-form over $Q$ which furthermore 
    is a quotient of $\varphi_{n}$, and
    the orthogonal quotient of $\varphi_{n}$ factors through it;
    \item if $\psi$ is dense on the right (resp.\ left) then $l$ 
    (resp.\ $r$) is the greatest left-sided (resp.\ right-sided)
    element below $n$;
    \item if $Q$ is unital then $\psi$ is dense on the right
    (resp.\ left) 
    if and only if $l$ 
    (resp.\ $r$) is the greatest left-sided (resp.\ right-sided)
    element below $n$.
    \end{enumerate}
\end{lemma}

\begin{proof}
    1. First we remark that the joins in 
    $\upsegment r$ are precisely the same as in $Q$, except for the 
    join of the empty set, which in $\upsegment r$ is $r$.
    Similarly, the joins 
    in $\upsegment l$ are those of $Q$ but with the empty join being 
    $l$. Hence, for $\psi$ to be a 2-form it suffices to verify that 
    it satisfies
    $\psi(r,y)=\psi(x,l)=0$ for all $y\in\upsegment l$ and 
    $x\in\upsegment r$. But we have $\psi(r,y)=\varphi_{n}(r,y)=0$ if 
    and only if $r\then y\le n$, which is true because $r$ is 
    right-sided: $r\then y\le r\le n$. Similarly, 
    $\psi(x,l)=0$ because $l$ is left-sided and $l\le n$,
    and we conclude that $\psi$ is a 2-form. Since it 
    is a quotient of $\varphi_{n}$, which is a 2-form over $Q$,
    $\psi$ is also a
    2-form over $Q$. Finally, the orthogonal quotient of 
    $\varphi_{n}$ is the least quotient of $\varphi_{n}$ and thus
    factors through 
    $\psi$.\\
    2. Now assume that $\psi$ is dense on the right, and let
    $a\le n$ be 
    a left-sided 
    element of $Q$. Then $a\vee l\le n$, and
    $1\then (a\vee l)\le n$, \ie,
    $\psi(1,a\vee l)=0$ (this makes sense because $a\vee 
    l\in\upsegment l$). Hence, since $\psi$ is 
    dense it follows that $a\vee l=l$, \ie, $a\le l$, 
    which shows that $l$ is the greatest left-sided
    element below 
    $n$. The situation with density on the left is similar.\\
    3. Assume that $Q$ is unital and that $l$ is the greatest
    left-sided element below $n$. Let $x\in\upsegment l$
    such that
    $1\perp x$, \ie, such that $1\then x\le n$.
    Then $1\then x\le l$ 
    because $1\then x$ is left-sided, and thus
    $x\le l$, \ie\ $x=l$, because 
    $Q$ is unital and thus $x\le 1\then x$.
    This shows that $\psi$ is 
    dense on the right. Density on the left is handled
    similarly.
    \qed
\end{proof}

\begin{theorem}
    Let $Q$ be a quantale, and
    $\varphi:L\times R\rightarrow 2$ a principal 2-form over $Q$
    with 
     generators $x\in L$ and $y\in R$. Then:
    \begin{enumerate}
    \item if $\varphi$ is dense on the right (resp.\ left) then 
    $\ann(y)$
    (resp.\ $\ann(x)$) is the greatest left-sided
    (resp.\ right-sided)
    element below $\orth(x,y)$;
    \item if $Q$ is unital then $\varphi$ is dense on the
    right (resp.\ left) 
    if and only if $\ann(y)$ 
    (resp.\ $\ann(x)$) is the greatest left-sided
    (resp.\ right-sided)
    element below $\orth(x,y)$.
    \end{enumerate}
\end{theorem}

\begin{proof}
    Let $n=\orth(x,y)$. From (\ref{thm:principalforms})
    it follows that $\varphi$ is a $Q$-orthoquotient of 
    $\varphi_{n}$, and
    (\ref{prop:densequotient}) implies that
    this quotient factors through
    $(f,g):\psi\rightarrow\varphi$, where the 2-form
    $\psi:\upsegment\ann(x)\times\upsegment\ann(y)\rightarrow 2$
    is as in (\ref{lemma:denseforms}), and both $f$
    and $g$ are 
    surjective and dense. Hence, by 
    (\ref{prop:pseudo}) we conclude that
    $\varphi$ is dense on the right if and only if $\psi$ is, and 
    similarly on the left. The result now follows from 
    (\ref{lemma:denseforms}).
    \qed
\end{proof}

\section{Involutive modules}\label{sec:invmods}

If $Q$ is an involutive quantale and we have a 2-form 
$\varphi:L\times R\rightarrow 2$ where both $L$ 
and $R$ are left $Q$-modules, it still makes sense to define
when it is that $\varphi$ is 
balanced, for the involution makes $L$ a 
right module: $xa=a^{*}x$. Hence, being balanced 
corresponds to the condition
$\form{a^{*}x}y=\form x{ay}$ for all $x\in L$, $y\in 
R$, and $a\in Q$, or, equivalently,
$\form{ax}{y}=\form x{a^{*}y}$.
We will not pursue this 
in general, but rather study the particular situation where $L=R$ 
and the 2-form is symmetric. Since in this situation we have
$\rorth{(-)}=\lorth{(-)}$, we shall 
write $\rorth{(-)}$ for both.

\begin{definition}
    Let $Q$ be an involutive quantale, $M$ a left $Q$-module, and 
    $\varphi$ a symmetric 2-form on $M$.
    The pair $(M,\varphi)$ (or just $M$, when no confusion may
    arise)
    is an \emph{involutive 
    (left) $Q$-module}
    if for all $a\in Q$ and 
    $x,y\in M$ we have
    \[\form{a^{*}x}y=\form x{ay}\;.\]
    A \emph{homomorphism} of involutive left $Q$-modules is a 
    homomorphism of left $Q$-modules $f$ such that $(f,f)$ is an
    orthomorphism:
    $\form {f(x)}{f(y)}=\form x y$.
    An \emph{involutive right module} is defined analogously by
    the condition
    \[\form{xa}y=\form x{ya^{*}}\;.\]
\end{definition}

The fact that we have restricted to symmetric 2-forms allows us to 
use the fact (\ref{prop:symvsinv}) that $\Q(\varphi)$ is an
involutive quantale:

\begin{proposition}
    Let $Q$ be an involutive quantale, and 
    $\varphi:M\times M\rightarrow 2$ a symmetric 2-form.
    There is a bijection between involutive left $Q$-module 
    structures on $(M,\varphi)$ and involution preserving
    homomorphisms
    from $Q$ to $\Q(\varphi)$.
\end{proposition}

\begin{proof}
    A quantale homomorphism 
    $h:Q\rightarrow\Q(\varphi)$ is the same as an action of $Q$ on 
    $\varphi$, with
    $h(a)=((-)a,a(-))$. Hence, $h$ preserves involution
    if and only if
    $((-)a^{*},a^{*}(-))=((-)a,a(-))^{*}=(a(-),(-)a)$, \ie,
    if and 
    only if $(-)a=a^{*}(-)$ for all $a\in Q$, \ie, if and only
    if $(M,\varphi)$ is an 
    involutive left $Q$-module. \qed
\end{proof}

From here and (\ref{exm:endomorphisms}) we see that in the case
when the
2-forms involved are faithful
the notion of involutive module corresponds precisely to that of
involutive 
representation
$Q\rightarrow\Q(S)$ of~\cite{MulvPell01,PellRosi97}. In other
words 
we have, as a corollary of (\ref{prop:equivalentformsoverQ}):

\begin{proposition}
    Let $Q$ be an involutive quantale, $M$ a left $Q$-module, 
    and $\varphi$ a symmetric 2-form on $M$. Then $(M,\varphi)$
    is an involutive left 
    $Q$-module if and only if
    $\rorth{(a^*x)}=a\backslash (\rorth x)$
    for all $a\in Q$ and $x\in M$.
    In the case when $\varphi$ is faithful this condition is 
    equivalent to
    $a^*x=\rorth{(a\backslash (\rorth x))}$.
\end{proposition}

All the previous definitions and results can be specialized to the 
case of involutive modules. We highlight just a few facts:

\begin{proposition}
    Let $Q$ be an involutive quantale, $m$ a left-sided element, and 
    $n$ a self-adjoint element such that $m\le n$. Then the left 
    $Q$-module $\upsegment m$ is involutive, with the symmetric 
    2-form being defined by $a\perp b\iff a^{*}\then b\le n$.
\end{proposition}

\begin{proof}
    Immediate consequence of (\ref{lemma:denseforms}),
    because $m^{*}$ is right-sided and $m^{*}\le n^{*}=n$. \qed
\end{proof}

\begin{proposition}
    Let $Q$ be an involutive quantale, $M$ an involutive left 
    $Q$-module, and $x\in M$. Then $\orth(x,x)$ is self-adjoint.
\end{proposition}

\begin{proof}
    For all $a\in Q$ we have:
    \[\begin{array}{l}a\le\orth(x,x)\iff x\perp ax\iff a^{*}x\perp x\iff
    x\perp a^{*}x\iff\\
    \iff a^{*}\le\orth(x,x)\;. \qed
    \end{array}\]
\end{proof}

\begin{proposition}\label{prop:upsegannx}
    Let $Q$ be an involutive quantale, and $M$ an involutive left 
    $Q$-module with a  generator $x\in M$. Then 
    $\upsegment\ann(x)$ is an involutive left $Q$-module with 
    the symmetric 2-form defined by $a\perp_{x} b\iff a^{*}\then 
    b\le\orth(x,x)$, and
    the map $\upsegment \ann(x)\rightarrow M$ defined by
    $a\mapsto ax$ is a surjective homomorphism of involutive left
    $Q$-modules.
\end{proposition}

\begin{proof}
    Let us just see that the map $a\mapsto ax$ is an orthomorphism:
    \[a\perp_{x} b\iff a^{*}\then b\le\orth(x,x)\iff
    x\perp a^{*}bx\iff
    ax\perp bx\;. \qed\]
\end{proof}

\begin{corollary}
    Let $Q$, $M$, and $x$ be as in the previous proposition.
    If the symmetric 
    2-form of $\upsegment\ann(x)$ is faithful then
    $\upsegment\ann(x)$ and $M$ are isomorphic as involutive left 
    $Q$-modules.
\end{corollary}

We conclude by pointing out a few immediate consequences of these
results in the 
case of quantales associated to C*-algebras.
Recall from~\S\ref{sec:background} that if $A$ is a unital 
C*-algebra then the set of all the closed linear subspaces
of $A$ is a 
unital involutive
quantale $\Max A$. We then have:

\begin{lemma}\label{lem:classification}
    Let $A$ be a unital C*-algebra, $m$ a maximal left-sided element 
    of $\Max A$, and $n=m\vee m^{*}$.
    Then the symmetric 2-form on $\upsegment m$ determined by $n$ is
    faithful.
\end{lemma}

\begin{proof}
    From C*-algebra theory we know that there is a unique pure state 
    $\varphi:A\rightarrow\CC$ whose kernel is $n$, and that the
    quotient
    $H=A/m$ is a Hilbert space with 
    inner product defined by
    $\langle a+m,b+m\rangle=\varphi(b^{*}a)$. Hence, two vectors
    $a+m,b+m\in H$ are orthogonal if and only if $b^{*}a\in n$,
    \ie,
    $a^{*}b\in n$, and
    thus the isomorphism $f:\upsegment m\rightarrow \proj(H)$
    is an 
    orthomorphism because, for all $c,d\in\upsegment m$,
    $f(c)$ is orthogonal to $f(d)$
    in the lattice 
    $\proj(H)$
    of closed linear subspaces of $H$ if and only if
    $c^{*}\then d\le 
    n$. Therefore the 2-form on $\upsegment m$ is faithful
    because the 
    2-form on $\proj(H)$ is. \qed
\end{proof}

\begin{theorem}
    Let $A$ be a unital C*-algebra, and $M$ an involutive left 
    $\Max A$-module with a  generator $x$.
    Assume also that $\ann(x)$ is a maximal left-sided element.
    Then $M$ is isomorphic as an involutive left $\Max A$-module to 
    $\upsegment\ann(x)$.
\end{theorem}

\begin{proof}
    Let $m=\ann(x)$. The topological left $A$-module structure of
    $A/m$ makes
    $\proj(A/m)$ a left $\Max A$-module
    (\cf\ \S\ref{sec:background}).
    Furthermore, $A/m$ is involutive 
    in the sense that $\langle ax,y\rangle=\langle x,a^{*}y\rangle$, 
    and from here it follows easily that $\proj(A/m)$
    is involutive as a
    $\Max A$-module. Hence, we have a surjective homomorphism
    $\upsegment m\rightarrow M$ of 
    involutive left $\Max A$-modules, which must be an
    isomorphism 
    because the 2-form on $\upsegment m$ is faithful. \qed
\end{proof}

Following the terminology of~\cite{MulvPell01}, let us
define a \emph{Hilbert representation} of $\Max A$ to be an 
involutive left $\Max A$-module isomorphic to one of the form
$\proj(H)$
determined by a representation 
of $A$ on $H$, in the manner of \S\ref{sec:background}.
Then we obtain:

\begin{corollary}\label{cor:classification}
    Let $A$ be a unital C*-algebra, and $M$ an involutive left 
    $\Max A$-module with a  generator $x$.
    Assume also that $\ann(x)$ is a maximal left-sided element.
    Then $M$ is a Hilbert representation determined by an 
    irreducible representation of $A$.
\end{corollary}

\begin{proof}
    This follows from the previous results and the fact that for a 
    maximal ideal $m$ the quotient $A/m$ defines an irreducible 
    representation. \qed
\end{proof}
    
The existence of a generator whose annihilator is a maximal 
left-sided element is, as was already mentioned in
(\ref{rem:MPpoints}), equivalent to the property known as
non-triviality in~\cite{MulvPell01}, and
the above corollary corresponds 
to one of the implications in~\cite[Th.\ 9.1]{MulvPell01}.
The main difference between the proof in~\cite{MulvPell01}
and what we have done above is that we have used 2-forms.
Also, we have focused 
less on the properties of those elements of $\Max A$ known as
``pure 
states'' and instead more on the annihilators of the generators
of principal 
$\Max A$-modules, \eg, formulating non-triviality directly in
terms of the 
annihilators.
In~\cite[Th.\ 9.1]{MulvPell01} it is further assumed that the 
module $M$ is an \emph{algebraically irreducible} representation, 
\ie, that $M$ is atomic as a sup-lattice and that each atom is a 
 generator (equivalently, $M$ is atomic and everywhere principal).
 Therefore our present formulation is
more general, even though it is not so in an essential way because
in the proof of~\cite[Th.\ 9.1]{MulvPell01} the extra
conditions are not 
used. In~\cite{MulvPell01} it is further conjectured that every 
algebraically irreducible representation of $\Max A$ is
non-trivial.

\end{document}